\documentstyle[12pt]{article}
\chardef\REST="16
\chardef\lcoR="70
\chardef\rcoR="71
\chardef\PreQ="34
\newcommand{\RESTR}{{\msas\REST}}
\newcommand{\restriction}{\RESTR}
\newcommand{\lcor}{\mathopen{\msas\lcoR}}
\newcommand{\rcor}{\mathclose{\msas\rcoR}}
\newcommand{\preQ}{\mathbin{\msas\PreQ}}
\chardef\tglE="45
\chardef\tgrE="44
\chardef\vtgL="43
\chardef\vtgR="42
\newcommand{\vtgl}{\mathbin{\msas\vtgL}}

\font\bfit=cmbxti10 scaled \magstep1
\newcommand{\xiieusm}{\font\teusm=eusm10 scaled \magstep1
\font\seusm=eusm8 \font\zeusm=eusm8}
\newcommand{\xieusm}{\font\teusm=eusm10 scaled \magstephalf
\font\seusm=eusm8 \font\zeusm=eusm8}
\newcommand{\xeusm}{\font\teusm=eusm10 
\font\seusm=eusm7 \font\zeusm=eusm7}

\newcommand{\skr}[1]{{\mathchoice {\hbox{\teusm{#1}}} 
{\hbox{\teusm{#1}}}{\hbox{\seusm{#1}}}{\hbox{\zeusm{#1}}} }}
\newcommand{\xiieurm}{\font\teurm=eurm10 scaled \magstep1
\font\seurm=eurm8 \font\zeurm=eurm6}

\newcommand{\xiimsb}{\font\tmsbm=msbm10 scaled \magstep1
\font\smsbm=msbm8 \font\zmsbm=msbm8}
\newcommand{\ximsb}{\font\tmsbm=msbm10 scaled \magstephalf
\font\smsbm=msbm8 \font\zmsbm=msbm8}
\newcommand{\xmsb}{\font\tmsbm=msbm10 
\font\smsbm=msbm7 \font\zmsbm=msbm7}

\newcommand{\BBB}[1]{{\mathchoice {\hbox{\tmsbm{#1}}} 
{\hbox{\tmsbm{#1}}}{\hbox{\smsbm{#1}}}{\hbox{\zmsbm{#1}}} }}
\newcommand{\xiimsa}{\font\tmsam=msam10 scaled \magstep1
\font\smsam=msam8 \font\zmsam=msam8}
\newcommand{\ximsa}{\font\tmsam=msam10 scaled \magstephalf
\font\smsam=msam8 \font\zmsam=msam8}

\newcommand{\msas}[1]{{\mathchoice {\hbox{\tmsam{#1}}} 
{\hbox{\tmsam{#1}}}{\hbox{\smsam{#1}}}{\hbox{\zmsam{#1}}} }}
\def\addto#1#2{
\ifx\zone\undefined\let\zone=#1\def#1{\zone#2}\else
\ifx\ztwo\undefined\let\ztwo=#1\def#1{\ztwo#2}\else
\ifx\zthree\undefined\let\zthree=#1\def#1{\zthree#2}\else
\fi\fi\fi
}
\addto\normalsize   {\xiieusm\xiimsb\xiimsa\xiieurm}
\addto\small        {\xieusm\ximsb\ximsa}
\addto\footnotesize {\xmsb\xeusm}
\newtheorem{theorem}             {Theorem}
\newtheorem{corollary}  [theorem]{Corollary}
\newtheorem{claim}      [theorem]{Claim}
\newtheorem{definition} [theorem]{Definition}
\newtheorem{lemma}      [theorem]{Lemma}
\newtheorem{proposition}[theorem]{Proposition}
\newtheorem{remark}
{{\bfit What\/ $\IST$ knows about standard sets ?\/}}
\newtheorem{que}   {{\bfit Question.}}         
\newtheorem{prF}   {Proof.}                    

\newcommand{\bcor}{\begin{corollary}}
\newcommand{\ecor}{\end{corollary}}
\newcommand{\bcl} {\begin{claim}}
\newcommand{\ecl} {\end{claim}}
\newcommand{\bdf} {\begin{definition}\rm}
\newcommand{\edf} {\end{definition}} 
\newcommand{\ble} {\begin{lemma}}
\newcommand{\ele} {\end{lemma}}
\newcommand{\bpro}{\begin{proposition}} 
\newcommand{\epro}{\end{proposition}} 
\newcommand{\baq} {\begin{acknow}\rm}
\newcommand{\eaq} {\end{acknow}} 
\newcommand{\bqe} {\begin{que}\rm}
\newcommand{\eqe} {\end{que}} 
\newcommand{\brem}{\begin{remark}\rm}
\newcommand{\erem}{\end{remark}} 
\newcommand{\bpf} {\begin{prF}\rm}
\newcommand{\epf} {\qeD\end{prF}} 
\newcommand{\bte} {\begin{theorem}}
\newcommand{\ete} {\end{theorem}}
\newcommand{\epF}[1]{\qeD\ {({\it#1\/}\hspace{0.5ex})}\end{prF}} 
\newcommand{\bde}{\begin{description}}
\newcommand{\ede}{\end{description}}
\newcommand{\ben}{\begin{enumerate}}
\newcommand{\een}{\end{enumerate}}
\newcommand{\bit}{\begin{itemize}}
\newcommand{\eit}{\end{itemize}}
\newcommand{\bce}{\begin{center}}
\newcommand{\ece}{\end{center}}
\newcommand{\bay}{\begin{array}}
\newcommand{\eay}{\end{array}}
\newcommand{\bqu}{\begin{quotation}}
\newcommand{\equ}{\end{quotation}}
\newcommand{\bbox}[1]{{\bf{#1}}}

\newcommand{\IST} {\bbox{IST}}
\newcommand{\ZFC} {\bbox{ZFC}}
\newcommand{\ZF}  {\bbox{ZF}}
\newcommand{\zfo} {\bbox{ZFGC}}

\newcommand{\ii}  {\bbox{i}}

\newcommand{\fT}  {\bbox{T}}
\newcommand{\rbox}[1]{{\rm{#1}}}
\newcommand{\st}   {{\tt st}\,}   %{\rbox{st}\,}
\newcommand{\for}  {\mathbin{\hspace{0.2ex}{\tt forc}\hspace{0.2ex}}}
\newcommand{\tru}[2]  {{\tt Truth\hspace{0.2ex}}^{#1}_{\in,#2}}
\newcommand{\dil}[1] {{\tt Def\hspace{0.2ex}}_{\in,<}(#1)}
\newcommand{\ie} {\hbox{{\it i.}\hspace{0.3ex}{\it e.}}}

\newcommand{\ax}  {{\upa x}}
\newcommand{\uS}  {{\upa{\hspace{0.2ex}\dS}}}
\newcommand{\sss}{\hspace{0.2ex}}
\newcommand{\aeq} {\mathbin{\sss\upa{\hspace{-0.6ex}=}\sss}}
\newcommand{\ain} {\mathbin{\sss\upa{\hspace{-0.6ex}\in}\sss}}
\newcommand{\pfn}{{\cP_{\rbox{fin}}}}
\newcommand{\al}  {\alpha} 
\newcommand{\da}  {\delta}

\newcommand{\ga}  {\gamma}

\newcommand{\la}  {\lambda} 
\newcommand{\sg}  {\sigma}
\newcommand{\Sg}  {\Sigma}

\newcommand{\vpi} {\varphi}
\newcommand{\om}  {\omega}

\chardef\sigmA ="1B

\newcommand{\bbb}{\hspace*{0.5pt}} 
\newcommand{\dvoj}[1]{\mathord{\bbb{\BBB #1}\bbb}}

\newcommand{\dI} {{\dvoj I}}
\newcommand{\dS} {{\dvoj S}}

\newcommand{\skrsp}{\hspace*{0.5pt}} 
\newcommand{\scri}[1]{\mathord{\skrsp\skr #1\skrsp}} 
\newcommand{\cA}{{\scri A}}
\newcommand{\cP}{{\scri P}}
\newcommand{\cU}{{\scri U}}
\newcommand{\cX}{{\scri X}}
\renewcommand{\sup}[2]
{\mathord{\kern 0.02em\vphantom{X}^{#2}\kern -0.08em #1}}
\newcommand{\sups}[2]
{\mathord{\kern 0.02em\vphantom{X}^{#2}\kern -0.15em #1}}
\newcommand{\supi}[2]
{\mathord{\kern 0.02em\vphantom{X}^{#2}\kern -0.27em #1}}
\newcommand{\upa}[1]{\sups{#1}\ast}

\newcommand{\upst}{^{\rbox{st}}}
\newcommand{\DS}{\displaystyle}
\newcommand{\Askip}{\hspace{0.1ex}}
\newcommand{\est}{{\DS\exists\upst}} 
\newcommand{\fst}{{\DS\forall\Askip\upst}}
\newcommand{\fstf}{{\DS\forall\Askip^{\rbox{stfin}}}}

\newcommand{\sq}{\subseteq}
\newcommand{\cj}{\mathbin{\hspace{1pt}\&\hspace{1pt}}}      

\newcommand{\sus} {{\exists\,}}
\newcommand{\kaz} {{\forall\,}}
\newcommand{\uu}  {{U\hspace*{0.3ex}}}

\newcommand{\imp} {\Longrightarrow} 
\newcommand{\eqv} {\Longleftrightarrow} 
 
\newcommand{\res} {\mathbin{\restriction}}
\newcommand{\we}  {{\mathbin{\kern 1.3pt ^\wedge}}}

\newcommand{\ti}  {\times}
\newcommand{\qeD} {\hfill{$\mtho\Box$}}
\newcommand{\qedd}{\hfill{$\mtho\dashv$}}
\newcommand{\dm}  {$$}

\newcommand{\mo}  {\models}

\newcommand{\iy}  {\infty}
\newcommand{\tgl}{\vtgl}     %\triangleleft}

\newcommand{\mem}[1]{\dd\in{#1}}
\newcommand{\ste}[1]{\dd{\tt st}{$\mtho\in$-{#1}}}
\newcommand{\ang} [1]{\langle #1\rangle}
\newcommand{\gon} [1]{\lcor #1\rcor}
\newcommand{\ans} [1]{\{\hspace{0.2mm}#1\hspace{0.2mm}\}}
\newcommand{\inst}  {{\cal L}_{\in,\st}}

\newcommand{\inle}  {{\cal L}_{\in,<}}

\newcommand{\dd}[2] {\hhsur\hbox{$\mtho{#1}$-#2}}

\newcommand{\vom}   {\vspace{1mm}}
\newcommand{\stkd}[2]{\ang{#1\hspace{1pt};\hspace{1pt}#2}}
\newcommand{\stkt}[3]{\ang{#1\hspace{1pt};\hspace{1pt}{#2}
,{#3}}}
\newcommand{\stkc}[4]{\ang{#1\hspace{1pt};\hspace{1pt}{#2}
,{#3},{#4}}}
\newcommand{\msur}  {\protect\hspace{-1\mathsurround}}

\newcommand{\hhsur} {\hspace{0.3\mathsurround}}
\newcommand{\mtho}  {\mathsurround=0mm}
\newcommand{\itla}  {\item\label}

\begin{document}

\normalsize

\title{What internal set theory knows about standard sets} 
\author{Vladimir Kanovei
\thanks
{\ 
\protect\begin{minipage}[t]{0.9\textwidth}
Moscow Transport Engineering Institute,\protect\\
{\tt kanovei@mech.math.msu.su} \ and \ 
{\tt kanovei@math.uni-wuppertal.de}.
\protect\end{minipage}\protect\vspace{2pt}
}
%\protect\\
\thanks
{\ 
Partially supported by a grant from DFG and visiting 
appointments from University of Wuppertal 
and Max Planck Institute at Bonn. 
}
\and
Michael Reeken
\thanks
{\ Bergische Universit\"at -- 
GHS Wuppertal. \ {\tt reeken@math.uni-wuppertal.de}}
}
\date{\today}
\maketitle

\begin{abstract}
We characterize those models of $\ZFC$ which are embeddable, 
as the class of all standard sets, in a model of 
internal set theory $\IST$.\\[2mm] 
\end{abstract}

%\newpage

\subsection*{Introduction}

In the early 60s Abraham Robinson demonstrated that nonstandard 
models of natural and real numbers could be used to 
interpret the basic notions of analysis in the spirit of 
mathematics of the 17-th and 18-th century, \ie\ including 
infinitesimal and infinitely large quantities. 

{\it Nonstandard analysis\/}, the field of mathematics which 
has been initiated by Robinson's idea, develops in two 
different versions. 

The {\it model theoretic\/} version, following the original 
approach, interprets ``nonstandard'' notions via nonstandard 
models in the $\ZFC$ universe.

On the other hand, the {\it axiomatic\/} version more radically 
postulates that the whole universe of sets 
(including all mathematical objects) is arranged in a 
``nonstandard'' way, so that it contains both the objects of 
conventional, ``standard'' mathematics, called standard,
and objects of different nature, called nonstandard. 
The latter type includes infinitesimal and infinitely large numbers, 
among other rather unusual objects. 

%There has been said quite a lot {\it pro\/} and {\it contra\/} 
%each of the two versions, but it seems that each of them 
Each of the two versions has 
its collective of adherents who use it as a working tool to 
develop nonstandard mathematics. 

The most of those who follow the axiomatic version use  
{\it internal set theory\/} $\IST$ of Nelson \cite{n77} as 
the basic set theory. This is a theory in the language $\inst$ 
(that is the language containing the membership $\in$ and the 
unary predicate of standardness $\st$ as the only atomic 
predicates) which includes all axioms of $\ZFC$ in the 
\mem language together with three principles that govern the 
interactions between {\it standard\/} (\ie\ those sets $x$ which 
satisfy $\st x$) and {\it nonstandard\/} objects in the set universe. 
(See below.) 

It is known that $\IST$ is an {\it equiconsistent\/} extension of 
$\ZFC.$ Moreover, $\IST$ is a {\it conservative\/} extension of 
$\ZFC,$ so that {\it an \mem sentence 
$\vpi$ is a theorem of\/ $\ZFC$ iff\/ $\vpi\upst$ is 
a theorem of\/ $\IST,$\/} where $\vpi\upst$ is the formal 
relativization of $\vpi$ to the class $\dS=\ans{x:\st x}$ of all 
standard sets.
This 
result, due to Nelson, is sometimes considered as a reason to view 
$\IST$ as a syntactical tool of getting $\ZFC$ theorems often in a 
more convenient way than traditional tools of $\ZFC$ (= the 
``standard'' mathematics) allow.

However working with $\IST$ one should be interested to know 
whether its axioms reflect some sort of mathematical reality. 
One could expect that the relations between $\ZFC$ and $\IST$ 
are similar to those between the real line and the complex plane,
so that each model of $\ZFC$ could be embedded, as the class of 
all standard sets, in a model of $\IST.$ However this is not the
case: we demonstrated in \cite{hyp1} that the least \mem model 
of $\ZFC$ is not embeddable in a model of $\IST.$ This observation 
leads us to the question:

\bit
\item
{\it 
which ``standard'' models (\ie\ transitive \mem models) of\/ $\ZFC$ 
can be embedded, as the class of all standard sets, in a model of\/  
$\IST$ ?}
\eit
Let $\zfo$ ($\ZF$ plus Global Choice) be the theory, in the 
language $\inle$ with the binary predicates $\in$ and $<$ as the 
only atomic predicates, containing all of $\ZFC$ (with the 
schemata of {\it Separation\/} and {\it Collection\/}, 
or {\it Replacement\/}, in $\inle$), together with 
the axiom saying that $<$ wellorders the universe in such a way 
that each initial segment is a set. 

Suppose that $M$ is a transitive set, ordered by a relation $<$ 
so that $\stkt M\in<$ models $\zfo.$ A set $T\sq M$ will be 
called {\it innocuous\/} for ${\stkt M\in<}$ if, for any sets 
$y\sq x\in M$ such that $y$ is definable~\footnote
{\ By formulas of $\inle$ (with parameters in $M$) 
plus $T$ as an extra predicate.} 
%\pagebreak[3]
in the structure $\stkc M\in<T,$ we 
have $y\in M.$ (Thus it is required that $T$ does not destroy 
Separation in $\stkt M\in<$ -- but it can destroy Collection.)

Note that every \dd\inle formula having sets in $M$ as parameters 
can be naturally considered as an element of $M.$ Let 
$\tru M{<}$ denote the set of all closed \dd\inle formulas true 
in $\stkt M\in<.$ 

\bte
\label m
Let\/ $M$ be a transitive\/ \mem model of\/ $\ZFC.$ Then the 
existence of a wellordering\/ $<$ of\/ $M,$ such that\/ 
$\stkt M\in<$ models\/ $\zfo$ and\/ $\tru{M}{<}$ is innocuous 
for\/ $\stkt M\in<,$ is necessary and sufficient for\/ $M$ to 
be embeddable, as the class of all standard sets, in a model 
of\/ $\IST$.~\footnote
{\ The \dd\IST embeddable transitive models of $\ZFC$ can be 
characterized in different terms. Suppose that $\stkt M\in<$ is a 
model of $\zfo.$ Let $\cX$ be a collection of subsets of $M.$ 
%Consider $\stkd{\stkt M\in<}\cX$ as a second order structure. 
Say that $\cX$ is {\it innocuous\/} for ${\stkt M\in<}$ if 
we have $y\in M$ whenever $y\sq x\in M$ and $y$ is definable in 
the second order structure $\stkd{\stkt M\in<}\cX.$
Then, a transitive model\/ $M\models \ZFC$ is embeddable, as 
the class of all standard sets, in a model of\/ $\IST,$ iff there 
is a wellordering\/ $<$ of\/ $M$ such that\/ $\stkt M\in<$ 
models\/ $\zfo$ and the family\/ $\cX$ of all sets\/ $X\sq M,$ 
definable in\/ $\stkt M\in<,$ is innocuous for\/ $\stkt M\in<$.

The equivalence of this characterization and the one given 
by the theorem can be easily verified directly without a 
reference to the \dd\IST embeddability.
}
\ete

This is the main result of the paper. 

The proof of the sufficiency is a modification of the original 
construction of an $\IST$ model by Nelson~\cite{n77}. The 
necessity is more interesting: it is somewhat surprising that 
$\IST$ ``knows'' that the standard universe is a model of $\zfo.$ 
On the other hand the involvment of the truth relation could 
be expected in view of the fact that $\IST$ provides a uniform 
truth definition for \mem formulas, see Theorem~\ref{tru} below.

\brem
$\,$\\[0ex]
The theorem answers the question in the title as follows: 
\bit
\item\msur
{\it $\IST$ ``knows'' about the standard universe that it can 
be wellordered by a relation\/ $<$ which respects the\/  
$\ZFC$ schemata of Separation and Collection, and moreover, the 
truth relation for the universe endowed by\/ $<$ does not destroy 
Separation.} 
\eit
This observation could perhaps lead to new insights in the 
philosophy of nonstandard mathematics. \hfill $\mtho\dashv$
\erem
 
It would be interesting to get similar results for other known 
nonstandard set theories, including those of Hrbacek~\cite{h78,mj} 
%(revised as $\HST$ in our survey~\cite{mj}) 
and Kawa\"\i~\cite{kaw83}.

\subsection{Internal set theory}
\label{ist}

Internal set theory $\IST$ is a theory in the language $\inst$ 
containing all axioms of $\ZFC$ (in the \mem language) and the 
following ``principles'':
\bde
\itemsep=1mm
\item{{\bfit Transfer\/}{\bf:}} \ $\sus x\:\Phi(x)\imp\est x\:\Phi(x)$ 
 \\ --- for any \mem formula $\Phi(x)$ with standard parameters;

\item{{\bfit Idealization\/}{\bf:}} \ 
$\fstf A\:\sus x\:\kaz a\in A\:\Phi(a,x)\eqv \sus x\:\fst a\:
\Phi(a,x)$ 
 \\ ---  for any \mem formula $\Phi(x)$ with arbitrary parameters;

\item{{\bfit Standardization\/}{\bf:}} \ 
$\fst X\:\est Y\:\fst x\:(x\in Y\eqv x\in X\cj \Phi(x))$ \\ ---  
for any \ste formula $\Phi(x)$ with arbitrary parameters.
\ede
The quantifiers $\est x$ and $\fst x$ have the obvious meaning 
(there exists a standard set $x$...). $\fstf A$ means: for any 
standard finite set $A$.

We shall systematically refer to different results in $\IST$ 
from \cite{n77,rms,hyp1}. In particular we shall use the following 
theorem of \cite{rms}.

\bte
\label{tru}
There is a\/ \ste formula\/ $\tau(x)$ such that, for any 
\mem formula\/ $\vpi(x_1,...,x_n),$ it is a theorem of\/ $\IST$ 
that\/ 
\dm
\fst x_1\:...\:
\fst x_n\:(\vpi\upst(x_1,...,x_n)\eqv
\tau(\gon{\vpi(x_1,...,x_n)}))\,.
\eqno{\Box}
\dm
\ete
Here $\gon\psi$ is the formula $\psi$ considered as a finite 
sequence of (coded) symbols of the \mem language and sets which 
occur in $\psi$ as parameters.
Thus the \mem truth in $\dS$ can be expressed by a single 
\ste formula in $\IST$.

\subsection{The necessity}
\label{nec}

Let us fix a transitive model $\dS$ of $\ZFC$ which is the 
standard part of a model $\dI$ of $\IST.$ To set things precisely, 
both $\dS$ and $\dI$ are sets in the $\ZFC$ universe, $\dS$ is a 
transitive \mem model of $\ZFC,$ so that 
${\in_\dS}={{\in}\res{\dS}},$ $\dS\sq\dI,$ $\dI$ is a model of 
$\IST,$ ${\in_\dS}={{\in_\dI}\res{\dS}},$ 
but of course ${\in_\dI}\not={{\in}\res{\dI}}$.

Our aim is to prove that there is an ordering $<$ of $\dS$ such 
that $\stkt \dS\in<$ models $\zfo$ and $\tru{\dS}{<}$ is 
innocuous for $\stkt \dS\in<$.

\subsubsection{The forcing and generic structures}
\label f

Let $\Sg$ be the class of all structures of the form 
$\sg=\stkd X<,$ where $X\in\dS$ is transitive and has the form 
$X=\dS_\al=V_\al\cap\dS$~\footnote
{\ $\msur V_\al$ is the \dd\al th level of the von Neumann set 
hierarchy.} 
for some ordinal $\al\in\dS,$ and ${<}\in\dS$ is a wellordering of 
$X.$ 

We say that $\sg'=\stkd{X'}{<'}$ {\it extends\/} 
$\sg=\stkd{X}<,$ symbolically $\sg\preQ\sg',$ if 
$X\sq X'$ and $<'$ is an end-extension of $<$. 

Define a relation ${\sg\for\Phi(x_1,...,x_n)},$ where 
$\sg={\stkd{X}<}\in\Sg$ while $\Phi$ is a \dd{\inle}formula 
and $x_1,...,x_n\in X,$ by induction on the complexity 
of $\Phi$.

\ben
\itemsep=1mm
\itla{f1}
If $\Phi$ is an elementary formula of $\inle,$ \ie\ 
$x<y,\msur$ $x=y,$ or $x\in y,$ 
%and $\sg={\stkd X<}\in\Sg,$ 
then $\sg\for\Phi$ iff $\Phi$ is true in $\sg.$~\footnote
{\ Here and sometimes below any $\sg=\stkd X<\in\Sg$ is 
understood as $\stkt X\in<$.}

\itla{f2}\msur
$\sg\for(\Phi\cj\Psi)$ iff $\sg\for\Phi$ and $\sg\for\Psi$.

\itla{f3}\msur
$\sg\for(\neg\,\Phi)$ iff there does not exist $\sg'\in\Sg$ 
extending $\sg$ such that $\sg'\for\Phi$.

\itla{f4}\msur
$\sg\for\sus x\,\Phi(x)$ iff there is $x\in X$ such that 
$\sg\for\Phi(x)$.
\een

For a \dd\inle formula $\Phi,$ a structure 
$\sg={\stkd{X}<}\in\Sg$ is called \dd\Phi{\it complete\/} 
iff, for any subformula $\Psi(x_1,...,x_n)$ of $\Phi$ and all 
$x_1,...,x_n\in X,$ we have $\sg\for\Psi(x_1,...,x_n)$ or 
$\sg\for\neg\,\Psi(x_1,...,x_n)$.

\bte
\label{comp}
{\rm($\IST$)} \ If\/ $\Phi$ is a closed\/ \dd\inle formula 
with sets in\/ $X$ as parameters, and\/ $\sg={\stkd{X}<}\in\Sg$ 
is\/ \dd\Phi complete, then\/ $\sg\for\Phi$ iff\/ $\sg\models\Phi$.
\ete
\bpf
By metamathematical induction on the complexity of $\Phi$.
\epf

\subsubsection{Increasing sequence of structures}
\label q

We shall define an increasing sequence of structures 
$\sg_\ga={\stkd{X_\ga}{<_\ga}}\in\Sg,\msur$ $\ga<\la,$ such that 
$\dS=\bigcup_{\ga<\la}X_\ga,$ hence the relation 
${<}=\bigcup_{\ga<\la}{<_\ga}$ wellorders $\dS.$ The structures 
$\sg_\ga$ will be rather ``complete'' (in the sense above); 
then $<$ will not destroy Replacement by an elementary chain 
argument. 

We face, however, a problem at limit steps: how to guarantee 
that the unions of $<_\ga$ still belong to 
$\dS.$ Now $\dI$ enters the reasoning. 
It occurs that the construction can be maintained in 
$\dI,$ so that, by the $\IST$ axiom of 
Standardization, the unions at limit steps will be still in $\dS$! 

%To determine the construction at nonlimit steps, f
Fix sets $D$ and 
$\tgl$ in $\dI$ such that the following 
holds in $\dI:$ $\dS\sq D$ and $\tgl$ is a (strict) wellordering 
of $D.$ Then $\tgl$ may not be a wellordering of $D$ from the 
point of view of the $\ZFC$ universe $V,$ but 
still $\tgl$ wellorders any set $S\in\dS$ in $V$ by Standardization 
and the fact that $\dS$ is a transitive set. 

We say that a structure $\sg\in\Sg$ is {\it totally complete\/} 
if it is \dd\Phi complete for any formula $\Phi$ of $\inle.$ 
The construction depends on the frequency of totally complete 
structures in $\Sg$.\vom

{\bfit Case 1\/}{\bf:} each $\sg\in\Sg$ can be extended to a 
totally complete $\sg'\in\Sg$.\vom

Define a sequence of structures 
$\sg_\ga={\stkd{X_\ga}{<_\ga}}\in\Sg\msur$ $(\ga<\la)$ such that 
$X_\da=\bigcup_{\ga<\da}X_\ga$ and 
${<}_\da=\bigcup_{\ga<\da}{<}_\ga$ for all limit ordinals 
$\da<\la,$ and $\sg_{\ga+1}$ is the \dd\tgl least totally 
complete structure in $\Sg$ which properly extends $\sg_\ga$.

Let $\la$ be the largest ordinal such that $\sg_\ga$ is defined 
(and belongs to $\Sg,$ hence to $\dS$) for all $\ga<\la\,;$ 
clearly $\la\leq$ ``the least ordinal not in $\dS$''.\qedd\vom

{\bfit Case 2\/}{\bf:} otherwise. \vom

Fix a recursive enumeration $\ans{\Phi_n:n\in\om}$ of all 
formulas of $\inle.$ A structure $\sg\in\Sg$ will be called 
\dd n{\it complete\/} if it is \dd{\Phi_k}complete for any 
$k\leq n$. 

We set $\la=\om$ in this case, pick a structure $\sg_0\in\Sg$ 
not extendable to a totally complete structure, and define a 
sequence of structures 
$\sg_n={\stkd{X_n}{<_n}}\in\Sg$ such that, for any $n\in\om,$ 
$\sg_{n+1}$ is the \dd\tgl least \dd ncomplete structure 
in $\Sg$ which properly extends $\sg_n$.\qedd\vom

In each of the two cases $\ang{\sg_\ga:\ga<\la}$ is a 
sequence of elements of $\dS.$ It can hardly be 
expected that the sequence is \mem definable in $\dS$ as the 
construction refers to notions which involve the 
\mem truth relation for $\dS.$ But the following holds:

\bpro
\label d
The sequence\/ $\ang{\sg_\ga:\ga<\la}$ is\/ \ste definable in\/ 
$\dI$. 
\epro
\bpf
Apply Theorem \ref{tru}.
\epf

\subsubsection{The order}

First of all we prove 

\ble
\label{sl}
$\la$ is a limit ordinal and\/ $\bigcup_{\ga<\la}X_\ga=\dS$. 
\ele 
\bpf
Recall that $\la=\om$ in Case 2. If $\la=\ga+1$ in Case 1 then, by 
the assumption of Case 1, we would be able to define $\sg_\la.$  
Hence $\la$ is a limit ordinal  and the relation 
${<}=\bigcup_{\ga<\la}{<}_\ga$ is a wellordering of $X$.  

Suppose that $X=\bigcup_{\ga<\la}X_\ga\not=\dS.$ Then $X\in\dS$  
as any of $X_\ga$ has the form $\dS_\al=V_\al\cap\dS$ for some 
$\al.$ Now $<$ belongs to $\dS$ by Standardization, being 
\ste defin\-able in $\dI$ by Proposition~\ref{d}. It follows that 
$\sg=\stkd{X}<\in\Sg.$ Moreover $\sg$ is totally complete. (As the 
limit of an increasing sequence of totally complete structures in 
Case 1, and by similar reasons in Case 2.) This immediately 
contradicts the choice of $\sg_0$ in Case 2, while, in Case 1, 
adds an extra term to the sequence, which contradicts the 
choice of $\la$.
\epf

It follows that ${<}=\bigcup_{\ga<\la}{<}_\ga$ is a wellordering 
of $\dS$. 

\bcor
\label<
$\stkt\dS\in<$ is a model of\/ $\zfo$.
\ecor
\bpf
To see that $\dS$ satisfies Separation in the language $\inle$ 
note that $<$ is \ste definable in $\dI$ by Proposition~\ref d 
and apply Standardization in $\dI.$ Now consider Collection. 
Suppose that $p,\,X\in\dS$ and $\Phi(x,y,p)$ is a \dd\inle 
form\-ula.
%containing some $p\in\dS$ as the only parameter. 
We have to find 
$Y\in\dS$ such that the following is true in $\dS$:
\dm
\kaz x\in X\;
[\,\sus y\:\Phi(x,y,p)\imp \sus y\in Y\:\Phi(x,y,p)\,]\,.
\dm
In both Case 1 and Case 2, there is $\ga<\la$ such that 
$p,\,X\in X_\ga$ and $\sg_\ga$ is 
\dd{(\sus y\:\Phi(x,y,p))}complete. Prove that $Y=X_\ga$ is 
as required.

Consider $x\in X,$ hence $\in X_\ga.$ Suppose that there is 
$y\in\dS$ such that $\Phi(x,y,p)$ holds in $\dS,$ and prove that 
such a set $y$ exists in $X_\ga$.

It follows from Lemma~\ref{sl} that $\dS$ is the union of an 
increasing chain of \dd{(\sus y\:\Phi(x,y,p))}complete structures. 
Therefore, by Theorem~\ref{comp} and an ordinary model-theoretic 
argument, $\stkt\dS\in<$ is an elementary extension of 
$\stkt{X_\ga}\in{{<_\ga}}$ with respect to the formula 
${\sus y\:\Phi(x,y,p)}$ and all its subformulas. This proves the 
existence of $y$ in $X_\ga$.
\epf

\subsubsection{The set of true formulas is innocuous}

Let $T=\tru{\dS}{<}$ be the set of all closed \dd\inle formulas 
(with sets in $\dS$ as parameters) true in the model 
$\stkt\dS\in<.$ The next lemma completes the proof of the 
necessity part in Theorem~\ref m.

\ble
\label i
$T$ is innocuous for\/ $\stkt\dS\in<$.
\ele
\bpf
It suffices to check that $T$ is \ste definable in $\dI.$ (Then 
the result follows by Standardization in $\dI$ as above.) 

Let $\Phi(p_1,...,p_k)$ be a closed \dd\inle formula with 
parameters $p_1,...,p_k\in\dS.$ Let $n$ be the number of 
$\Phi(x_1,...,x_k)$ (see Case 2 in Subsection~\ref q). Take the 
least $\ga<\la$ such that $p_1,...,p_k\in X_\ga$ and, in Case 2, 
$\ga\geq n.$ Arguing as in the proof of Corollary~\ref<, we 
conclude that $\sg_\ga$ is an elementary substructure of 
$\stkt\dS\in<$ with respect to $\Phi,$ in particular 
$\Phi(p_1,...,p_k)$ is either true or false simultaneously in 
both $\sg_\ga$ and $\stkt\dS\in<.$ It 
remains to recall that the sequence of structures $\sg_\ga$ is 
\ste definable in $\dI$ by Proposition~\ref d.
\epf

\subsection{The sufficiency}
\label{suf}

This section proves the sufficiency part in Theorem~\ref m. 
We start with a transitive set $\dS$ and a wellordering $<$ 
of $\dS$ such that $\stkt\dS\in<\models\zfo,$ and suppose 
that the set $\fT=\tru{\dS}{<}$ of all closed \dd\inle formulas 
(with parameters in $\dS$), true in $\stkt\dS\in<,$ 
is innocuous for $\stkt\dS\in<.$ The aim is to embed $\dS,$ as the 
class of all standard sets, in a model $\dI$ of $\IST$.

\subsubsection{The ultrafilter}
\label u

To obtain $\dI$ we shall use the construction of an {\it adequate\/} 
ultrapower of Nelson~\cite{n77}, modified by Kanovei~\cite{rms}. 

Let $\dil\dS$ denote the collection of all sets $X\sq\dS$ 
definable in $\stkt\dS\in<$ by a formula of $\inle$ containing sets 
in $\dS$ as parameters.

Let $I=\pfn(\dS)=\ans{i\sq\dS:i\,\hbox{ is finite}}.$ This is a 
proper class in $\dS.$ 
%We shall use $I$ as the index set to define the ultrapower. 
Let $\cA$ be the algebra of all sets $X\sq I$ which belong to 
$\dil\dS$.

\bpro
\label U
There exists an ultrafilter\/ $U\sq\cA$ satisfying
\ben
\def\theenumi{{\rm(\Alph{enumi})}}
\def\labelenumi{\theenumi}
\itla A
if\/ $a\in\dS$ then the set\/ $\ans{i\in I:a\in i}$ belongs to\/ 
$U\;;$

\itla B
if\/ $P\sq\dS\times I,\msur$ $P\in\dil\dS,$ then the following set 
is in\/ $\dil\dS:$
\dm
\ans{x\in\dS:\hbox{\rm the cross-section 
$P_x=\ans{i:\ang{x,i}\in P}$ belongs to $U$}}\;;
\dm

\itla C
there is a set\/ $\cU\sq\dS,$ definable in the structure\/ 
$\stkc\dS\in<\fT,$ such that\/ $U=\ans{\cU_x:x\in\dS},$ where\/ 
$\cU_x=\ans{i\in I:\ang{x,i}\in\cU}$ for all $x$.
\een
\epro
\bpf
{\it Step $0.$\/} Let $U_0$ be the collection of all sets of the 
form 
\pagebreak[1]
\dm
I_{a_1...a_m}=\ans{i\in I:a_1,...,a_m\in i},
\hspace{3mm}\hbox{where}\hspace{3mm}a_1,..,a_m\in\dS\,.
\dm
The family $U_0$ obviously satisfies FIP 
({\it the finite intersection property\/}).\vom\qedd

{\it Step $n+1.$\/} Suppose that a FIP family $U_n$ of subsets 
of $I$ has been constructed. Denote by $\chi_n(x,i)$ the \dd nth 
formula in a recursive enumeration, fixed beforehand, of all 
\dd\inle formulas with exactly two free variables. 

We define $U_{n+1}=U_n\cup\ans{B_x:x\in\dS},$ where 
$B_x$ is equal to the set 
$A_x=\ans{i\in\dI:\stkt\dS\in<\models\chi_n(x,i)}$ whenever the 
family $U_n\cup\ans{B_y:y<x}\cup A_x$ still satisfies FIP, and 
$B_x=I\setminus A_x$ otherwise.\vom\qedd

Clearly $U=\bigcup_n U_n$ is as required. \ref C follows from the 
fact that the whole construction can be carried out in 
$\stkc\dS\in<\fT$.
\epf

Let us fix such an ultrafilter $U\sq\cA.$ 

Let $\uu i\:\Phi(i)$ mean: 
``the set $\ans{i\in I:\stkt\dS\in<\mo\Phi(i)}$ belongs to $U$''. 
(The quantifier: there exist \dd Umany.) Then, by the choice of 
$U,$ we have 
$\uu i\:(a\in i)$ for any $a\in\dS,$ and, given a relation 
$P(i,...)$ in $\dil\dS,$ the relation $\uu i\:P(i,...)$ belongs 
to $\dil\dS$ as well. 

\subsubsection{The model}
\label M

For $r\geq 1,$ we let $I^r=I\ti...\ti I$ ($r$ times $I$), and 
\dm
F_r=\ans{f\in\dil\dS:f\,\hbox{ maps }\,I^r\,\hbox{ to }\,\dS}\,.
\dm
Let separately $I^0=\ans{0}$ and 
$F_0=\ans{\ans{\ang{0,x}}:x\in\dS}$. 

We finally put $F_\iy=\bigcup_{r\in\om} F_r,$ and, for 
$f\in F_\iy,$ let $r(f)$ be the only $r$ such that $f\in F_r$.

Suppose that $f\in F_\iy,\msur$ $q\geq r=r(f),$ and 
$\ii=\ang{i_1,...,i_r,...,i_q}\in I^q.$ Then we set 
$f[\ii]=f(\ang{i_1,...,i_r}).$ In particular $f[\ii]=f(\ii)$ 
whenever $r=q.$ Separately we put $f[\ii]=x$ for any $\ii$ 
whenever $f=\ans{\ang{0,x}}\in\dS_0$.

Let $f,\,g\in F_\iy$ and $r=\max\ans{r(f),r(g)}.$ Define 
\dm
f\aeq g
\hspace{4mm}\hbox{iff}\hspace{4mm}
\uu i_r\;\uu i_{r-1}\;...\;\uu i_1\;(f[\ii]=g[\ii]),
\dm
(where $\ii$ denotes $\ang{i_1,...,i_r}$), 
and define $f\ain g$ similarly. (Note the order of quantifiers.) 
The following is a routine statement.

\bpro
\label r
\ $\aeq$ is an equivalence relation on\/ $F_\iy.$ The relation\/ 
$\ain$ on\/ $F_\iy$ is\/ \dd\aeq invariant in each of the two 
arguments.
\qeD
\epro

Define $[f]=\ans{g\in F_\iy:f\aeq g}.$ Let 
$\uS=\ans{[f]:f\in F_\iy}$ (the quotient). 
For $[f],\,[g]\in\uS,$ define $[f]\ain[g]$ iff $f\ain g.$ 
(This is independent of the choice of representatives by the 
proposition.) 

For any $x\in\dS,$ define $\ax=[\ans{\ang{0,x}}],$ the image of 
$x$ in $\uS$. 

We finally define $\st[f]$ iff $[f]=\ax$ for some $x\in\dS.$ 

\bte
\label S
$\stkt\uS\ain\st$ is a model of\/ $\IST.$ The map\/ 
$x\longmapsto\ax$ is a 1--1\/ \mem embedding of\/ $\dS$ onto the 
class of all standard elements of $\uS$.
\ete

The theorem immediately implies the sufficiency part in 
Theorem~\ref m.

\bpf
We begin with an appropriate formalism. Let $\Phi(f_1,...,f_m)$ 
be an \mem formula with functions $f_1,...,f_m\in F$ as 
parameters. Put $r(\Phi)=\max\ans{r(f_1),...,r(f_m)}.$ If 
$r\leq q$ and $\ii\in I^q$ then let $\Phi[\ii]$ denote the 
formula $\Phi(f_1[\ii],...,f_m[\ii])$ (an \mem formula with 
parameters in $\dS$). 
Let finally $[\Phi]$ denote $\Phi([f_1],...,[f_m]),$ which is an 
\mem formula with parameters in $\uS$.

\bpro
\label l
{\rm(\L o\v s)} \ 
Let\/ $\Phi=\Phi(f_1,...,f_m)$ be an \mem formula with functions 
$f_1,...,f_m\in F$ as parameters, and\/ $r=r(\Phi).$ Then
\dm
[\Phi]\,\hbox{ holds in }\,\uS
\hspace{6mm}\hbox{iff}\hspace{6mm}
\uu i_r\;\uu i_{r-1}\;...\;\uu i_1\;
(\Phi[\ii]\,\hbox{ holds in }\,\dS)\,.
\dm
\epro  
\bpf
($\ii$ denotes $\ang{i_1,...,i_r}$ in the displayed line.) The only 
detail one needs to note is that, since the index set $I$ is a 
proper class in $\dS,$ we need the global choice to carry out the 
ordinary argument. This is why $\dS$ needs to be a model of 
$\zfo,$ not merely $\ZFC$.
\epF{Proposition \ref l}

Using functions in $F_0,$ we immediately conclude that the map 
$x\longmapsto \ax$ is an \mem elementary 1--1 embedding of $\dS$ 
onto the class of all standard sets in $\uS,$ which implies both 
Transfer and all of $\ZFC$ axioms in $\stkt\uS\ain\st.$ It remains 
to check Idealization and Standardization.\vom

{\it Idealization.\/} Let $\Phi(a,x)$ be an \mem formula with two 
free variables, $a$ and $x,$ and some functions in $F$ as 
parameters. We have to demonstrate
\dm
\fstf A\:\sus x\:\kaz a\in A\:[\Phi](a,x)\;\imp\;
\sus x\:\fst a\:[\Phi](a,x)
\eqno{(\dag)}
\dm
in $\uS.$ (It is known that the implication $\Longleftarrow$ here 
is a corollary of other axioms of $\IST$.) 
The left--hand side of $(\dag)$ implies, by Proposition~\ref l, 
\pagebreak[0]
\dm
\forall_{\rm finite}\:A\sq\dS\;\uu i_r\;\uu i_{r_1}\;...\;\uu i_1\;
\sus x\;\kaz a\in A\;\Phi[\ang{i_1,...,i_r}](a,x)
\dm
in $\dS,$ where $r=r(\Phi).$ To simplify the formula note that the 
leftmost quantifier is a quantifier over $I$ and define a function 
$\al\in F_{r+1}$ by $\al(i_1,...,i_r,i)=i.$ The last displayed 
formula takes the form 
\dm
\kaz i\in I\;\uu i_r\;\uu i_{r_1}\;...\;\uu i_1\;
(\sus x\;\kaz a\in \al\;\Phi)[\ang{i_1,...,i_r,i}](a,x)\,,
\dm
which implies $\sus x\:\kaz a\in [\al]\:[\Phi](a,x)$ in $\uS$ by 
Proposition~\ref l. Now, by the definition of the predicate $\st$ in 
$\uS,$ it suffices to check that $\ax\ain[\al]$ in $\uS$ for any 
$x\in\dS.$ This is equivalent to 
$\uu i\;\uu i_r\;...\;\uu i_1\;(x\in i),$ 
which holds by the choice of $U$.\vom

{\it Standardization.\/} Recall that $U$ is definable in the 
structure $\stkc\dS\in<\fT$ by \ref C of Proposition~\ref U. 
Therefore the model $\stkt\uS\ain\st$ is definable in 
$\stkc\dS\in<\fT$ as well. Thus we have only to check that, given 
$x\in\dS,$ any set $y\sq x,$ which is definable in 
$\stkc\dS\in<\fT,$ 
%(with parameters from $\dS$ allowed) 
belongs to $\dS.$ But this follows from the fact that $\fT$ is 
innocuous for $\stkt\dS\in<$.

\epF{Theorems \ref S and \ref m}

\end{document}